\documentclass{article}
\usepackage{curves,epic,epsfig,makeidx,amsmath,amsfonts,latexsym}
\newtheorem{theorem}{Theorem}[section]
\newtheorem{prop}{Proposition}[section]
\newtheorem{lemma}{Lemma}[section]
\begin{document}

\title{\large\textbf{{ASYMPTOTIC SHAPE IN A CONTINUUM\linebreak GROWTH MODEL}}}
\author{Maria Deijfen \thanks{Postal address: Department of Mathematics, Stockholm
University, SE-106 91 Stockholm, Sweden. E-mail:
mia@matematik.su.se.}}
\date{November 2002}

\maketitle

\thispagestyle{empty}

\begin{abstract}
\noindent A continuum growth model is introduced. The state at
time $t$, $S_t$, is a subset of $\mathbb{R}^d$ and consists of a
connected union of randomly sized Euclidean balls, which emerge
from outbursts at their center points. An outburst occurs
somewhere in $S_t$ after an exponentially distributed time with
expected value $|S_t|^{-1}$ and the location of the outburst is
uniformly distributed over $S_t$. The main result is that if the
distribution of the radii of the outburst balls has bounded
support, then $S_t$ grows linearly and $S_t/t$ has a non-random
shape as $t\rightarrow \infty$. Due to rotation invariance the
asymptotic shape must be a Euclidean ball.

\vspace{1cm}


\vspace{0.5cm}

\hspace*{4.8cm} Secondary 82B43.
\end{abstract}

\section{Introduction}

\noindent There is a variety of random growth models defined in a
discrete space; see for instance Eden (1961), Williams and
Bjerknes (1972), Bramson and Griffeath (1981) and Lee and Cowan
(1994). A particular one is the Richardson model, introduced in
Richardson (1973). It describes a Markov process whose state at
time $t$, $S_t$, is a subset of $\mathbb{Z}^d$. Each site in
$\mathbb{Z}^d$ is in either of two states, denoted 0 and 1, and
$S_t$ consists of the sites which are in state 1 at time $t$. A
site in state 0 is transferred to state 1 at a rate proportional
to the number of nearest neighbors in state 1, and once in state 1
it never returns to state 0. Thus, if sites in state 1 are thought
of as infected sites and sites in state 0 as uninfected this
dynamics defines a pure growth model. The main result, first
proved in Richardson (1973), states that, if $S_0$ consists of a
single site, then $S_t/t$ has a non-random shape as $t\rightarrow
\infty$. Generalizations of Richardson's result can be found in
Cox and Durrett (1981), Kesten (1986) and Boivin (1990). Apart
from the fact that the asymptotic shape is convex and compact, not
much is known about its qualitative features.\medskip

\noindent In this paper we introduce a growth model defined in a
continuous space, that is, the state at time $t$, still denoted by
$S_t$, is a subset of $\mathbb{R}^d$ instead of $\mathbb{Z}^d$.
The process should be thought of as describing the spread of some
kind of infection in a continuous medium and just as in the
Richardson model, the set $S_t$ specifies the region infected at
time $t$. The growth takes place by way of outbursts in the
infected region. Given the development of the infection up to time
$t$, the time until an outburst occurs somewhere in $S_t$ is
exponentially distributed with parameter $|S_t|$ and the location
of the outburst is uniformly distributed over $S_t$ (the notation
$|\cdot|$ will, throughout the paper, be used to denote both
Lebesgue measure and Euclidean norm, but the meaning will always
be obvious from the context). When an outburst occurs at an
infected point it causes a ball of stochastic radius around the
outburst point to be infected and the total infected region is
enlarged by the amount of this ball that was not previously
infected. Consequently the infected region is a connected union of
Euclidean balls. The radii of these balls are assumed to be
i.i.d.\ with distribution $F$ and mean $\gamma$.\medskip

\noindent The main result in this paper is a shape theorem for the
continuum model. An essential advantage of the continuum model as
compared to the Richardson model is that it possesses rotational
invariance, which forces the asymptotic shape to be a Euclidean
ball. To formulate the theorem, let $B(x,r)$ denote a ball with
radius $r$ around the point $x\in \mathbb{R}^d$ and let $S_t$
denote the infected region at time $t$ in the $d$-dimensional
continuum model starting at time zero from a ball with radius
$\gamma$ around the origin.

\begin{theorem}[Shape theorem] \label{th:aform}
Assume that $F$ has bounded support and let $S_0=B(0,\gamma)$.
Then, for any dimension $d$, there is a real number $\mu>0$ such
that, for any $\varepsilon$ with $0<\varepsilon<{\mu}^{-1}$,
almost surely
$$
(1-\varepsilon)B(0,{\mu}^{-1})\subset \frac{S_t}{t}
\subset(1+\varepsilon)B(0,{\mu}^{-1})
$$
\noindent for all sufficiently large $t$.
\end{theorem}

\noindent Another example of a continuum model with a Euclidean
ball as asymptotic shape is described in Howard and Newman
(1997).\medskip

\noindent The condition that the support of $F$ is bounded can
probably be weakened but some assumption on the tail is certainly
necessary to ensure that $\mu>0$, i.e. to ensure that the growth
is not faster than linear.\medskip

\noindent A natural question is how the asymptotic shape is
affected if the region infected at time zero is chosen to be
something other than a ball with radius $\gamma$ around the
origin. The answer is given in the following theorem.

\begin{theorem}\label{th:oberS0}
The conclusion of the shape theorem holds for any bounded initial
set $S_0$ with strictly positive Lebesgue measure. Also, the
number $\mu$ is independent of the initial set.
\end{theorem}

\noindent The rest of the paper is organized as follows. We start
by describing the model more thoroughly in Section 2. In Section 3
we define some important quantities and prove a technical lemma.
Section 4 contains results concerning the growth in a fixed
direction. To obtain these results we formulate and prove a
sharpened version of Liggett's (1985) subadditive ergodic theorem
(Theorem \ref{th:liggmod}). Section 5 is devoted to the proof of
the shape theorem and the proof of Theorem \ref{th:oberS0}.

\section{Description of the model}

\noindent In this section we construct the model more formally by
defining a Markov process whose state at time $t$, $S_{t}$, is a
subset of $\mathbb{R}^d$. The process may be thought of as
describing the spread of an infection (with no recoveries) or the
growth of a germ colony in a continuous medium. Points in $S_{t}$
will be referred to as infected.\medskip

\noindent To define the model, let $N$ be a Poisson process on
${\mathbb{R}}^{d+1}$. The extra dimension represents the time
dimension and the points of $N$ are hence denoted $(X_k,T_k)$,
where $X_k\in {\mathbb{R}}^{d}$ and the last coordinate $T_k$
gives the location on the time axis. To each point in the Poisson
process we associate a random radius $R_k$. The variables
$\{R_k\}$ are assumed to be i.i.d.\ with mean $\gamma$. At time
zero a ball with radius $\gamma$ around the origin, denoted by
$B_0$, is infected. The idea now is to follow the cylinder
$B_0\times \mathbb{R}$ upwards along the time axis until a point
in the Poisson process is found. An outburst then takes place at
this point generating a new infection ball $B_1$, whose size is
given by the radius associated with this particular Poisson point.
The new infected region is given by $B_0\cup B_1$. Scanning within
the cylinder $(B_0\cup B_1)\times \mathbb{R}$ further upwards
along the time axis we eventually hit a new Poisson point
representing a new outburst and corresponding enlargement of the
infected region. And so on.\medskip

\noindent To make this description more formal, for $S\subset
{\mathbb{R}}^{d}$, let $N_{S\times \mathbb{R}}$ denote the
restriction of $N$ to $S\times {\mathbb{R}}$. The growth of the
infected area takes place at time points $\{T_n\}$ by aid of
outbursts with radii $\{R_n\}$ at points $\{X_n\}$ obtained from
the following recursion:

\begin{itemize}
\item[1.] Let $X_0=0$, $T_0=0$, $R_0=\gamma$ and define $B_n=\{y\in {\mathbb{R}}^{d};
\hspace{0.1cm}|X_n-y|\leq R_n\}$, that is, $B_n$ is a ball in
${\mathbb{R}}^{d}$ with radius $R_n$ centered at $X_n$.
\item[2.] Given $\{X_i; \hspace{0.1cm} i\leq n\}$,
$\{T_i; \hspace{0.1cm} i\leq n\}$ and $\{R_i; \hspace{0.1cm} i\leq
n\}$, the time $T_{n+1}$ is defined as
$$
T_{n+1}=\inf_{k} \{T_k;\hspace{0.1cm} T_k>T_n \textrm{ and }
(X_k,T_k)\in N_{\cup_{i=0}^{n}{B_i}\times {\mathbb{R}}}\}
$$
\noindent and $X_{n+1}$ is the (a.s.\ unique) point in
${\mathbb{R}}^{d}$ such that $(X_{n+1},T_{n+1})\in
N_{\cup_{i=0}^{n}{B_i}\times {\mathbb{R}}}$. The radius $R_{n+1}$
is given by the radius associated with $(X_{n+1},T_{n+1})$.
\end{itemize}

\noindent From the sequence $\{X_n\}$ a new sequence $\{S_{(n)}\}$
is constructed by defining $S_{(n)}=\cup_{i=0}^{n}{B_i}$. The
infected region at time $t$ is now given by
$$
S_t=S_{(n)} \quad \textrm{for } t\in [T_n,T_{n+1}).
$$
\noindent Let us introduce the notation $\Delta_n=T_n-T_{n-1}$,
$n\geq 1$, for the successive times between the outbursts. By
construction of the model and properties of the Poisson process
$$
\Delta_{n+1}|{\mathcal{F}}_{n}  \sim \textrm{Exp}(|S_{(n)}|),
$$
\noindent where ${\mathcal{F}}_n=\sigma (X_0,\ldots
,X_{n},T_0,\ldots ,T_n, R_0,\ldots ,R_n)$. Given $S_t$, the
memoryless property of the exponential distribution implies that
the time until an outburst occurs somewhere in $S_t$ is
exponentially distributed with parameter $|S_t|$. The location of
the outburst, $X_n$ -- where $n$ is such that $t\in [T_{n-1},T_n)$
-- is uniformly distributed over $S_t$ and the new infected region
is given by $S_t\cup B_n$. Furthermore, the model is
Markovian.\medskip

\noindent To guarantee that the model is defined for all $t$ one
detail remains to be checked: We have to make sure that the
sequence $\{T_n\}$ does not have a finite limit point
$T_{\infty}$, since this would cause problems defining $S_t$ for
$t>T_{\infty}$. The following proposition is what we need:

\begin{prop}
Assume that the radii distribution has finite moment of order $d$.
Then, almost surely, $T_n \rightarrow \infty$ as $n\rightarrow
\infty$.
\end{prop}

\noindent \emph{Proof:} Let $\{E_k\}$ be independent, $E_k\sim
\textrm{Exp}(k)$. It is left to the reader to show that it
suffices to prove that $\sum_{k=1}^{\infty}E_k=\infty$ with
probability 1. To establish this, introduce
$\tilde{E}_k=E_k-\textrm{E}[E_k]=E_k-1/k$. Using the fact that
$\sum_{k=1}^{\infty}\textrm{E}[{\tilde{E}_k}^2]=
\sum_{k=1}^{\infty}1/k^2<\infty$, Kolmogorov's three-series
theorem implies that $\sum_{k=1}^{\infty} \tilde{E}_k$ converges
almost surely. Thus
$$
\sum_{k=1}^{\infty}E_k=\sum_{k=1}^{\infty}\tilde{E}_k+\sum_{k=1}^
{\infty}\frac{1}{k}=\infty,
$$
\noindent since $\sum_{k=1}^{\infty}1/k=\infty$.$\hfill\Box$

\section{Preliminaries}

\noindent In this section we introduce some preliminary notation
and results needed to prove the shape theorem. \medskip

\noindent To begin with, let $T(x)$ denote the time when the point
$x$ is infected, i.e.
$$
T(x)=\inf \{t; \hspace{0.1cm} x\in S_t\}.
$$
\noindent Our first result is a lemma bounding the time it takes
for the infection to reach the point $x$. Here, for
$z\in\mathbb{R}$, $\lceil z\rceil$ is the smallest integer larger
than $z$.

\begin{lemma} \label{lemma:txyb}
For any $x\in \mathbb{R}^d$ there exist i.i.d.\ exponential
variables $\{E_k\}$ with parameter $\lambda=\lambda(d)$ such that
$$
T(x)\leq \sum_{k=1}^{2\lceil|x|\gamma^{-1}\rceil}E_k.
$$
\end{lemma}

\noindent \emph{Proof:} (We give the proof for $d=2$. The case
$d\geq 3$ is analogous.)\medskip

\noindent Fix $x\in \mathbb{R}^2$ and assume, without loss of
generality, that $x$ is located on the $x$-axis, that is, assume
$x=(x',0)$.\medskip

\noindent To begin with, fix $c\in \mathbb{R}$ and write $B_k$ for
the ball of radius $c$ around the point $(k\cdot \gamma /2,0)$,
see Figure 1. If $c$ is sufficiently small, say $c\leq \gamma
/10$, then $B_1$ is contained in $B(0,\gamma)$, that is, $B_1$ is
infected at time zero. Let $E_1$ denote the time from time zero
until an outburst with radius at least $\gamma$ occurs in $B_1$.
Clearly $B_2$ is infected by the time such an outburst has
occurred, that is, $B_2\subset S_{E_1}$. Now, the outburst points
whose radius exceeds $\gamma$ constitutes a Poisson process with
intensity $p=P(R_k\geq \gamma)$ and thus $E_1\sim
\textrm{Exp}(p\pi c^2)$, where $\pi c^2$ is the area of
$B_1$.\medskip

\noindent The idea of how to continue should be clear: Let $E_0=0$
and define $E_k$, $k\geq 1$, recursively as the time until an
outburst with radius larger than $\gamma$ occurs in $B_k$ counting
from time $E_0+\ldots +E_{k-1}$. Since $B_k\subset S_{E_1+\ldots
+E_{k-1}}$, we have that $E_k\sim \textrm{Exp}(p\pi c^2)$. When an
outburst has occurred in $B_{2\lceil |x|\gamma^{-1}\rceil}$ it is
clear that $x$ must be infected, i.e. $x\in S_{E_1+\ldots
+E_{2\lceil |x|\gamma^{-1}\rceil}}$. Hence
$$
T(x)\leq \sum_{k=1}^{2\lceil |x|\gamma^{-1}\rceil} E_k
$$
\noindent and the lemma is proved. $\hfill \Box$\medskip

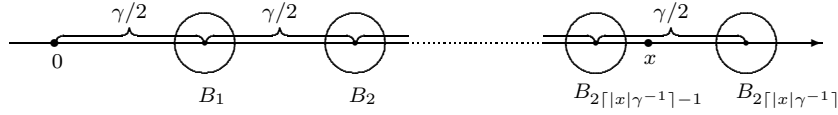
\begin{figure}
\centering \setlength{\unitlength}{2cm} \linethickness{0.7pt}
\begin{picture}(5.1,1.2)(-2,-0.6)

\thinlines \put(-2.3,0){\line(1,0){2.65}}
\dottedline{0.035}(0.35,0)(1.25,0)
\put(1.25,0){\vector(1,0){1.85}}

\bigcircle[0]{0.4} \scaleput(-1,0){\bigcircle[0]{0.4}}
%\scaleput(-2,0){\bigcircle[0]{0.4}}
\scaleput(1.6,0){\bigcircle[0]{0.4}}
\scaleput(2.6,0){\bigcircle[0]{0.4}}

\put(-2,0){\circle*{0.05}}
\put(-1.98,-0.12){\makebox(0,0){\footnotesize $0$}}
\put(-1,0){\circle*{0.03}} \put(0,0){\circle*{0.03}}
\put(1.6,0){\circle*{0.03}} \put(2.6,0){\circle*{0.03}}
\put(1.95,0){\circle*{0.05}}
\put(1.96,-0.11){\makebox(0,0){\footnotesize $x$}}

\curve[0](-1.99,0,-1.9754,0.0354,-1.94,0.05)
\put(-1.94,0.05){\line(1,0){0.39}}
\curve[0](-1.55,0.05,-1.5146,0.06464,-1.5,0.1)
\curve[0](-1.5,0.1,-1.4854,0.06464,-1.45,0.05)
\put(-1.45,0.05){\line(1,0){0.4}}
\curve[0](-1.055,0.05,-1.0196,0.0354,-1.005,0)
\put(-1.5,0.2){\makebox(0,0){\footnotesize $\gamma/2$}}

\curve[0](-0.99,0,-0.9754,0.0354,-0.94,0.05)
\put(-0.94,0.05){\line(1,0){0.39}}
\curve[0](-0.55,0.05,-0.5146,0.06464,-0.5,0.1)
\curve[0](-0.5,0.1,-0.4854,0.06464,-0.45,0.05)
\put(-0.45,0.05){\line(1,0){0.4}}
\curve[0](-0.055,0.05,-0.0196,0.0354,-0.005,0)
\put(-0.5,0.2){\makebox(0,0){\footnotesize $\gamma/2$}}

\curve[0](1.61,0,1.6246,0.0354,1.66,0.05)
\put(1.66,0.05){\line(1,0){0.39}}
\curve[0](2.05,0.05,2.0854,0.06464,2.1,0.1)
\curve[0](2.1,0.1,2.1146,0.06464,2.15,0.05)
\put(2.15,0.05){\line(1,0){0.4}}
\curve[0](2.545,0.05,2.5804,0.0354,2.595,0)
\put(2.1,0.2){\makebox(0,0){\footnotesize $\gamma/2$}}

\curve[0](0.01,0,0.0246,0.0354,0.06,0.05)
\put(0.06,0.05){\line(1,0){0.29}}

\put(1.25,0.05){\line(1,0){0.3}}
\curve[0](1.545,0.05,1.5804,0.0354,1.595,0)

\put(-0.95,-0.35){\makebox(0,0){\footnotesize $B_1$}}
\put(0.05,-0.35){\makebox(0,0){\footnotesize $B_2$}}
\put(1.88,-0.35){\makebox(0,0){\footnotesize $B_{2\lceil
|x|\gamma^{-1}\rceil -1}$}}
\put(2.88,-0.35){\makebox(0,0){\footnotesize $B_{2\lceil
|x|\gamma^{-1}\rceil}$}}

\end{picture}
\caption{A chain of small balls located $\gamma /2$ units apart on
the line segment joining the origin and $x$ is constructed.}
\end{figure}

\noindent \textbf{Remark 3.1} It follows from the proof of the
lemma that the bound for $T(x)$ is valid also for the time until a
small cube around $x$ is infected. By construction of the
variables $\{E_k\}$, at time $\sum_{k=1}^{2\lceil |x|\gamma^{-1}
\rceil}E_k$ an outburst with radius at least $\gamma$ has occurred
somewhere in a ball with radius $c/2$, $c\leq \gamma/10$, centered
at a point within distance $\gamma/2$ from $x$. Thus, a cube
centered at $x$ will be contained in the infected area at time
$\sum_{k=1}^{2\lceil |x|\gamma^{-1}\rceil}E_k$ if the side length
of the cube is chosen small enough. This observation will be
useful in proving Lemma \ref{lemma:Stilb}.\medskip

\noindent We now introduce some auxiliary quantities that will be
of use later. Let $S_t^{(x,s)}$, $t\geq s$, denote the set of
points that can be reached from $x$ within time $t$ if a new
process is started at $x$ at time $s$. That is, at time $s$ all
infection except a ball with radius $\gamma$ around $x$ is erased.
The infection then evolves in time according to the same rules as
for the original process, using the same $(d+1)$-dimensional
Poisson process. This gives rise to a new process, emanating from
$x$, whose state at time $t$, $t\geq s$, is given by
$S_t^{(x,s)}$. Now let
$$
\tilde{T}(x)=\inf\{t; \hspace{0.1cm} B(x,\gamma)\subset S_t\}
$$
\noindent and
$$
\tilde{T}(x,y)=\inf\left\{t; \hspace{0.1cm} B(y,\gamma)\subset
S_{\tilde{T}(x)+t}^{(x,\tilde{T}(x))}\right\}.
$$
\noindent In words, $\tilde{T}(x)$ is the time when the entire
ball with radius $\gamma$ around $x$ is infected and
$\tilde{T}(x,y)$ is the time it takes for the infection to invade
the entire $\gamma$-ball around $y$ if a new process is started at
$x$ at time $\tilde{T}(x)$. Note that the quantity
$\tilde{T}(x,y)$ is independent of $\tilde{T}(x)$ and has the same
distribution as $\tilde{T}(y-x)$. Furthermore, it is clear that,
if a point is contained in the region infected at time
$\tilde{T}(x)+t$ in the process started at $x$ at time
$\tilde{T}(x)$, then it is also contained in the region infected
at time $\tilde{T}(x)+t$ in the original process, i.e.
$$
S_{\tilde{T}(x)+t}^{(x,\tilde{T}(x))} \subset S_{\tilde{T}(x)+t}
$$
\noindent and hence

\begin{equation} \label{eq:allmsubad}
\tilde{T}(y)\leq \tilde{T}(x)+\tilde{T}(x,y).
\end{equation}

\noindent We close this section by anticipating that in what
follows it will be convenient to have special notation for the
quantity $\tilde{T}(mx,nx)$. Thus, let
$$
\tilde{T}(mx,nx)=\tilde{T}_{m,n}(x).
$$
\noindent Since $\tilde{T}(nx)=\tilde{T}_{0,n}(x)$, the
subadditivity property (\ref{eq:allmsubad}) translates into

\begin{equation}\label{eq:subadd}
\tilde{T}_{0,n}(x)\leq \tilde{T}_{0,m}(x)+\tilde{T}_{m,n}(x).
\end{equation}

\section{Growth in a fixed direction}

\noindent The proof of the shape theorem basically consists of two
parts:

\begin{itemize}
\item[1.]Show that $S_t$ grows linearly in each fixed direction
and that the asymptotic speed of the growth in each direction is
an almost sure constant. By the rotational invariance of
${\mathbb{R}}^{d}$ and the model, this constant must be the same
for all directions.
\item[2.]Show that the linear growth of $S_t$ holds for all directions
simultaneously.
\end{itemize}

\noindent This section is devoted to the first part. The first
task is to prove the following result.

\begin{prop} \label{prop:gvtt} For each $x\in \mathbb{R}^d$, we have
\begin{itemize}
\item[\rm{(a)}] $\mu(x):=\lim_{n\rightarrow \infty}{\rm{E}}[\tilde{T}(nx)]/n=
\inf_{n\geq 1} {\rm{E}}[\tilde{T}(nx)]/n;$
\item[\rm{(b)}] $\lim_{n\rightarrow \infty}
\tilde{T}(nx)/n=\mu(x)$ a.s.
\end{itemize}
\end{prop}

\noindent \textbf{Remark 4.1} We employ the convention that limits
over $n$ are taken over the positive integers, while limits over
$t$, which occur later on in the paper, are taken over all
positive reals.\medskip

\noindent To prove Proposition \ref{prop:gvtt} we will invoke the
following theorem by Liggett (1985).

\begin{theorem}[Liggett's subadditive ergodic theorem] \label{th:Liggett}
Let $\{X_{m,n}\}$ be a collection of random variables indexed by
integers satisfying $0\leq m <n$. Suppose $\{X_{m,n}\}$ has the
following properties:

\begin{itemize}
\item[\rm{(i)}] $X_{0,n}\leq X_{0,m}+X_{m,n}$.
\item[\rm{(ii)}] For each $n$, {\rm{E}}$|X_{0,n}|<\infty$ and {\rm{E}}$[X_{0,n}]\geq
cn$ for some constant $c>-\infty$.
\item[\rm{(iii)}] The distribution of $\{X_{m,m+k}; \hspace{0.1cm} k\geq 1\}$
does not depend on $m$.
\item[\rm{(iv)}] For each $k\geq 1$, $\{X_{nk,(n+1)k}; \hspace{0.1cm} n\geq 0\}$ is a
stationary sequence.
\end{itemize}

\noindent Then

\begin{itemize}
\item[\rm{(a)}] $\eta :=\lim_{n\rightarrow \infty}
{\rm{E}}[X_{0,n}]/n=\inf_{n\geq 1} {\rm{E}}[X_{0,n}]/n$.
\item[\rm{(b)}] The limit $X:=\lim_{n\rightarrow \infty}
X_{0,n}/n$ exists a.s.
\item[\rm{(c)}]${\rm{E}}[X]=\eta$.
\end{itemize}

\noindent Furthermore, if the stationary processes in
$\mathrm{(iv)}$ are ergodic, then

\begin{itemize}
\item[\rm{(d)}] $X=\eta$ a.s.
\end{itemize}

\end{theorem}

\noindent A brief outline of the structure of the proof can be
found below. For more detail we refer the reader to Liggett
(1985).\medskip

\noindent We wish to apply this result to the variables
$\{\tilde{T}_{m,n}(x)\}$. It turns out, however, that condition
(iv) fails for the sequence $\{\tilde{T}_{nk,(n+1)k}(x);
\hspace{0.1cm} n\geq 0\}$. Luckily the assumption (iv) can be
relaxed without weakening the conclusions:

\begin{theorem}\label{th:liggmod}
Let $\{X_{m,n}\}$ be a collection of random variables satisfying
{\rm{(i)-(iii)}} of Theorem \ref{th:Liggett}. Furthermore, suppose
that

\begin{itemize}
\item[\rm{(iv')}] $\limsup_{n\rightarrow \infty}X_{0,nk}/n\leq
{\rm{E}}[X_{0,k}]$ for each $k$.
\end{itemize}

\noindent Then {\rm{(a)-(d)}} of Theorem 4.1 hold for
$\{X_{m,n}\}$.
\end{theorem}

\noindent \textbf{Remark 4.2} A trivial modification of the proof
of Theorem \ref{th:liggmod} yields that, if we confine ourselves
with the conclusions (a)-(c), then (iv') can be replaced by the
still weaker assumption

\begin{itemize}
\item[(iv'')] $\textrm{E}\left[\limsup_{n\rightarrow
\infty}X_{0,nk}/n\right]\leq \textrm{E}[X_{0,k}]$ for each $k$.
\end{itemize}

\noindent To prove Theorem \ref{th:liggmod} requires some
knowledge of the structure of the proof of Theorem
\ref{th:Liggett}. This knowledge is provided in the following
brief sketch.\medskip

\noindent Write
$$
\bar{X}=\limsup_{n\rightarrow \infty} \frac{X_{0,n}}{n}
$$
\noindent and
$$
\b{$X$}=\liminf_{n\rightarrow \infty} \frac{X_{0,n}}{n}.
$$
\noindent The proof of Theorem \ref{th:Liggett} is broken up into
three steps:

\begin{itemize}
\item[(L1)] Prove that $\eta=\lim_{n\rightarrow \infty} \textrm{E}[X_{0,n}]/n= \inf_{n\geq 1}
\textrm{E}[X_{0,n}]/n$.
\item[(L2)] Prove that $\textrm{E}[\bar{X}]\leq \eta$, and if the stationary
processes in (iv) are ergodic, then $\bar{X}\leq \eta$ almost
surely.
\item[(L3)] Prove that ${\rm{E}}[\b{$X$}]\geq \eta$.
\end{itemize}

\noindent From (L2) and (L3) it follows that
$\textrm{E}[\b{$X$}]\geq \textrm{E}[\bar{X}]$. This implies that
$\b{$X$}$ and $\bar{X}$ are equal, since trivially $\b{$X$}\leq
\bar{X}$. Hence, once (L2) and (L3) are accomplished it is clear
that $X:=\lim_{n\rightarrow \infty} X_{0,n}/n$ exists with
probability 1. It also follows from (L2) and (L3) that
$\textrm{E}[X]=\eta$ and by (L1) that $\eta<\infty$. Furthermore,
if $\bar{X}\leq \eta$ -- which, according to (L2), for example is
the case if the sequences $\{X_{nk,(n+1)k}; \hspace{0.1cm} n\geq
0\}$ are ergodic -- then we can deduce that $X=\eta$ almost
surely.\medskip

\noindent For (L1)-(L3) we refer the reader to Liggett (1985). The
essential task for us is to identify the parts of the proof that
make use of the assumption (iv). Since this assumption is to be
replaced by (iv'), these are the parts that have to be modified.
The following table shows how the assumptions are used in the
different steps.\medskip

\begin{center}
\begin{tabular}{c|c}
Step & Assumptions used \\ \hline (L1) & (i), (iii) \\ (L2) &
(i)-(iv)
\\ (L3) & (i)-(iii)
\end{tabular}
\end{center}\medskip

\noindent Since the proofs of (L1) and (L3) do not use (iv), these
statements hold also if (iv) is dropped. In proving (L2) though,
Liggett uses the assumption (iv) and thus a modification of
Liggett's proof is necessary to establish that (L2) remains true
when (iv) is replaced by (iv'). The modification is described in
the following proof of Theorem \ref{th:liggmod}.\bigskip

\noindent \emph{Proof of Theorem \ref{th:liggmod}:} It suffices to
show that $\bar{X}\leq \eta$ almost surely, that is, it suffices
to show that

\begin{equation} \label{eq:Lm6}
\limsup_{n\rightarrow \infty} \frac{X_{0,n}}{n}\leq \eta \quad
\textrm{a.s.}
\end{equation}

\noindent To achieve this, fix $\delta >0$ and choose $k$ large so
that $\textrm{E}[X_{0,k}]/k\leq \eta+\delta$. We will show that
for all $j$,

\begin{equation} \label{eq:Lm3}
\limsup_{n\rightarrow \infty} \frac{X_{0,nk+j}}{nk+j}\leq
\frac{\textrm{E}[X_{0,k}]}{k} \quad \textrm{a.s.}
\end{equation}

\noindent which yields
$$
\limsup_{n\rightarrow \infty} \frac{X_{0,n}}{n}\leq
\frac{\textrm{E}[X_{0,k}]}{k} \quad \textrm{a.s.}
$$
\noindent Hence, once (\ref{eq:Lm3}) has been established it
follows from the choice of $k$ that
$$
\limsup_{n\rightarrow \infty} \frac{X_{0,n}}{n}\leq \eta+\delta
$$
\noindent and, since $\delta >0$ was arbitrary, this implies
(\ref{eq:Lm6}). To prove (\ref{eq:Lm3}), fix $j$ and use
subadditivity to get

\begin{equation} \label{eq:Lm2}
\frac{X_{0,nk+j}}{nk+j}\leq \frac{n}{nk+j} \cdot
\frac{X_{0,nk}}{n}+ \frac{n}{nk+j} \cdot \frac{X_{nk,nk+j}}{n}.
\end{equation}

\noindent By (iii), the distribution of $X_{nk,nk+j}$ depends only
on $j$ and, by (ii), the first moment is finite. Thus,
$$
\sum_{n=1}^{\infty} P(X_{nk,nk+j}>n\varepsilon)<\infty
$$
\noindent for all $\varepsilon >0$ and by the Borel-Cantelli lemma
this implies that

\begin{equation} \label{eq:Lm1}
\lim_{n\rightarrow \infty} \frac{X_{nk,nk+j}}{n}=0 \quad
\textrm{a.s.}
\end{equation}

\noindent Using (iv'), (\ref{eq:Lm3}) follows from (\ref{eq:Lm2})
and (\ref{eq:Lm1}).$\hfill \Box$ \bigskip

\noindent We are now in a position to prove Proposition
\ref{prop:gvtt}.\bigskip

\noindent \emph{Proof of Proposition \ref{prop:gvtt}:} Fix $x\in
\mathbb{R}^d$. Since $\{\tilde{T}_{m,n}(x)\}$ satisfies the
assumptions (i)-(iii) of Theorem \ref{th:Liggett}, the proposition
follows from Theorem \ref{th:liggmod} if we can show that (iv')
holds for $\{\tilde{T}_{m,n}(x)\}$, that is, if we can show that

\begin{equation}\label{eq:(iv')fTt}
\limsup_{n\rightarrow \infty} \frac{\tilde{T}_{0,nk}(x)}{n}\leq
\textrm{E}[T_{0,k}] \quad \textrm{for all } k.
\end{equation}

\noindent To do this it is necessary to introduce an auxiliary
sequence $\{\tilde{T}_{(i-1)k,ik}'(x);\newline \hspace{0.1cm}
i\geq 1\}$ defined recursively as follows:

\begin{itemize}
\item[] Let $\tilde{T}_{0,k}'(x)=\tilde{T}_{0,k}(x)$. For
$i\geq 2$, given $\{\tilde{T}_{(l-1)k,lk}'(x); \hspace{0.1cm}
l\leq i-1\}$, define
$$
\tilde{T}_{(i-1)k,ik}'(x)=\inf\left\{t; \hspace{0.1cm}
B(ikx,\gamma)\subset S_{\Phi_{i-1}^{k}+
t}^{((i-1)kx,\hspace{0.1cm} \Phi_{i-1}^{k})}\right\},
$$
\noindent where $\Phi_{i-1}^{k}=\sum_{l=1}^{i-1}
\tilde{T}_{(l-1)k,lk}'(x)$.
\end{itemize}

\noindent Remember that $S_{t}^{(x,s)}$ is the area infected at
time $t$, $t\geq s$, in a process started at time $s$ emanating
from a ball with radius $\gamma$ around $x$. Thus,
$\tilde{T}_{(i-1)k,ik}'(x)$ is the time when the ball with radius
$\gamma$ around $ikx$ is infected in a process started at
$(i-1)kx$ at time $\Phi_{i-1}^{k}$. Some thought reveals that the
variables $\{\tilde{T}_{(i-1)k,ik}'(x)\}$ are i.i.d.\ with
expected value $E[\tilde{T}_{0,k}(x)]$. Hence, by the strong law
of large numbers,

\begin{equation} \label{eq:Lm5}
\frac{1}{n}\sum_{i=1}^{n}\tilde{T}_{(i-1)k,ik}'(x) \rightarrow
\textrm{E}[\tilde{T}_{0,k}(x)] \quad \textrm{as } n\rightarrow
\infty.
\end{equation}

\noindent Furthermore, it is readily seen that
$$
\tilde{T}_{0,nk}(x)\leq \sum_{i=1}^{n} \tilde{T}_{(i-1)k,ik}'(x).
$$
\noindent Dividing this inequality by $n$ and using (\ref{eq:Lm5})
we obtain (\ref{eq:(iv')fTt}). Theorem \ref{th:liggmod} now gives
that all the conclusions (a)-(d) of Theorem \ref{th:Liggett} hold
for $\{\tilde{T}_{m,n}(x)\}$. Part (a) of the proposition follows
from \ref{th:Liggett}(a) and part (b) follows from Theorem
\ref{th:Liggett}(b) and (d).$\hfill \Box$ \bigskip

\noindent Our next task is to show that $\mu(x)$ is nonzero so
that the growth is indeed linear.

\begin{prop} \label{prop:0muinf}
If $F$ has bounded support, then $0<\mu(x)<\infty$ for each $x\in
\mathbb{R}^d$, $x\neq 0$.
\end{prop}

\noindent \emph{Proof:} Fix $x\in \mathbb{R}^d$. That
$\mu(x)<\infty$ follows from Proposition \ref{prop:gvtt}(a) so it
remains to show that $\mu(x)>0$. Since $F$ has bounded support,
there is a real number $r_{max}$ such that with probability 1 the
radius of an outburst ball does not exceed $r_{max}$. Using this,
it follows from a minor modification of Lemma 3.1 in Penrose
(2001)\footnote{or see the appendix on p. 47.} that there are
constants $c_1,c_2\in (0,\infty)$ such that

\begin{equation} \label{eq:PTxb}
P(T(x)\leq c_1|x|)\leq 2\cdot 3^{-c_2|x|}.
\end{equation}

\noindent That $\mu(x)$ is nonzero follows from (\ref{eq:PTxb})
and the fact that

\begin{equation} \label{eq:ETnxmu}
\limsup_{n\rightarrow \infty} \frac{\textrm{E}[T(nx)]}{n}\leq
\lim_{n\rightarrow \infty}
\frac{\textrm{E}[\tilde{T}(nx)]}{n}=\mu(x)
\end{equation}

\noindent by arguing as follows: Substituting $x$ by $nx$ in
(\ref{eq:PTxb}) yields
$$
P\left(\frac{T(nx)}{n}\leq c_1|x|\right)\leq 2\cdot 3^{-c_2n|x|},
$$
\noindent which implies that
$$
\frac{\textrm{E}[T(nx)]}{n}\geq c_1|x|\left(1-2\cdot
3^{-c_2n|x|}\right).
$$
\noindent Using (\ref{eq:ETnxmu}) this gives that $\mu(x)\geq
c_1|x|$. $\hfill \Box$\bigskip

\noindent The next step is to prove that the discrete limit in
Proposition \ref{prop:gvtt}(b) can be replaced by a continuous
one.

\begin{prop} \label{prop:kontpar}
For each $x\in \mathbb{R}^d$ we have $\lim_{t\rightarrow \infty}
\tilde{T}(tx)/t=\mu(x)$, where the limit is taken along $t\in
\mathbb{R}^+$.
\end{prop}

\noindent \emph{Proof:} Fix $x\in \mathbb{R}^d$. We start by
showing that

\begin{equation} \label{eq:tildeTncx}
\lim_{n\rightarrow \infty} \frac{\tilde{T}(nqx)}{nq}= \mu(x) \quad
\textrm{for all } q\in \mathbb{Q},
\end{equation}

\noindent that is, moving away from the origin in direction $x$
using steps of arbitrary rational length yields the same limit. To
this end, introduce the notation
$$
\overline{T}_{t}=\frac{\tilde{T}(tx)}{t}
$$
\noindent Now, since $q\in \mathbb{Q}$ we have $q=k/m$ for some
integers $k$ and $m$. The sequence $\{\overline{T}_{nk}\}$ is a
subsequence of $\{\overline{T}_n\}$ and, hence, by Proposition
\ref{prop:gvtt}(b),
$$
\lim_{n\rightarrow \infty}\overline{T}_{nk}=\lim_{n\rightarrow
\infty}\overline{T}_{n}= \mu(x) \quad \textrm{a.s.}
$$
\noindent However, $\{\overline{T}_{nk}\}$ is also a subsequence
of $\{\overline{T}_{nq}\}$ -- obtained by considering only those
points where $n$ is a multiple of $m$ -- and since the sequence
$\{\overline{T}_{nq}\}$ does indeed have a limit (this follows
from Proposition \ref{prop:gvtt}(b) applied to the point $qx$)
this implies that
$$
\lim_{n\rightarrow \infty}\overline{T}_{nq}=\lim_{n\rightarrow
\infty}\overline{T}_{nk}= \mu(x) \quad \textrm{a.s.}
$$
\noindent To complete the proof we use (\ref{eq:tildeTncx}) to
show that

\begin{equation} \label{eq:lstildeTtxmmuli}
\limsup_{t\rightarrow \infty}\frac{\tilde{T}(tx)}{t}\leq
\mu(x)\leq \liminf_{t\rightarrow \infty}\frac{\tilde{T}(tx)}{t}.
\end{equation}

\noindent To establish (\ref{eq:lstildeTtxmmuli}) we will need a
bound for the time from the moment when the $\gamma$-ball around
an arbitrary point on the line segment between $nqx$ and $(n+1)qx$
is infected until the $\gamma$-balls of all points on the line
segment are infected. To obtain such a bound, let $l_{nqx}$ denote
the line segment between $nqx$ and $(n+1)qx$ and write
$\tilde{T}(l_{nqx})$ for the time when all points on $l_{nqx}$ has
their $\gamma$-balls infected, i.e.
$$
\tilde{T}(l_{nqx})=\inf\{t; \hspace{0.1cm} B(z,\gamma)\subset S_t
\textrm{ for all } z\in l_{nqx}\}.
$$

\begin{figure}
\centering \setlength{\unitlength}{2.9cm} \linethickness{0.7pt}
\begin{picture}(2.6,2.2)(-1.3,-1.1)

\bigcircle[0]{2} \bigcircle[0]{1.44}
\curve[0](0.36,0.6235,0.5,0.866) \curve[0](0.6235,0.36,0.866,0.5)
\curve[0](0.72,0,1,0) \curve[0](0.6235,-0.36,0.866,-0.5)
\curve[0](0.36,-0.6235,0.5,-0.866) \curve[0](0,0.72,0,1)
\curve[0](0,-0.72,0,-1) \curve[0](-0.36,0.6235,-0.5,0.866)
\curve[0](-0.6235,0.36,-0.866,0.5) \curve[0](-0.72,0,-1,0)
\curve[0](-0.6235,-0.36,-0.866,-0.5)
\curve[0](-0.36,-0.6235,-0.5,-0.866)

\thinlines \put(-1.3,0){\line(1,0){2.6}}

\curve[0](0,0,0.5,0.866) \curve[0](0,0,0.5,-0.866)
\curve[0](0,0,0.866,0.5) \curve[0](0,0,0.866,-0.5)
\curve[0](0,0,0,1) \curve[0](0,0,0,-1) \curve[0](0,0,-0.5,0.866)
\curve[0](0,0,-0.5,-0.866) \curve[0](0,0,-0.866,0.5)
\curve[0](0,0,-0.866,-0.5)

\put(-0.6,0.6){\makebox(0,0){\footnotesize $\text{A}_1$}}
\put(-0.2,0.83){\makebox(0,0){\footnotesize $\text{A}_2$}}
\put(0.25,0.83){\makebox(0,0){\footnotesize $\text{A}_3$}}
\put(0.65,0.6){\makebox(0,0){\footnotesize $\text{A}_4$}}
\put(0.85,0.21){\makebox(0,0){\footnotesize $\text{A}_5$}}
\put(0.85,-0.21){\makebox(0,0){\footnotesize $\text{A}_6$}}
\put(0.65,-0.6){\makebox(0,0){\footnotesize $\text{A}_7$}}
\put(0.25,-0.83){\makebox(0,0){\footnotesize $\text{A}_8$}}
\put(-0.2,-0.83){\makebox(0,0){\footnotesize $\text{A}_9$}}
\put(-0.6,-0.6){\makebox(0,0){\footnotesize $\text{A}_{10}$}}
\put(-0.82,-0.21){\makebox(0,0){\footnotesize $\text{A}_{11}$}}
\put(-0.83,0.21){\makebox(0,0){\footnotesize $\text{A}_{12}$}}

\put(-0.22,0){\circle*{0.05}}
\put(-0.35,-0.08){\makebox(0,0){\footnotesize $nqx$}}

\put(0.25,0){\circle*{0.05}}
\put(0.45,-0.08){\makebox(0,0){\footnotesize $(n+1)qx$}}

\end{picture}
\caption{The front zone $F$ of $B(z_0,\gamma)$ divided into pieces
$A_1,\ldots ,A_{12}$.}
\end{figure}
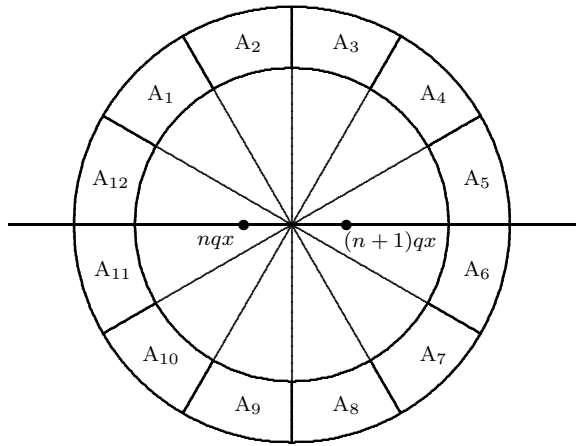

\noindent Assume that at time $t_0$ there is a point $z_0\in
l_{nqx}$ such that $B(z_0,\gamma)\subset S_{t_0}$. For small
$q\in\mathbb{Q}$ we will derive an upper bound for
$\tilde{T}(l_{nqx})-t_0$ expressed as the maximum of a number of
i.i.d.\ exponential random variables. In two dimensions such a
bound is easily obtained using a geometric construction displayed
in Figure 2. Namely, let $F$ be the outer ring of width
$\gamma/10$ in $B(z_0,\gamma)$, that is,
$F=B(z_0,\gamma)\backslash B(z_0,9\gamma/10)$. Divide $F$ into
twelve disjoint pieces $A_1,\ldots ,A_{12}$ of equal area as shown
in Figure 2 and let $E_k$ ($k=1,\ldots ,12$) be the time from time
$t_0$ until an outburst with radius at least $\gamma$ occurs in
$A_k$. By construction, the areas $A_1,\ldots ,A_{12}$ are all
infected at time $t_0$, implying that $E_k\sim
\textrm{Exp}(p|A_k|)$. Furthermore, since the $A_k$ are disjoint,
the variables $\{E_k\}$ are independent. Now, if $q$ is small,
then all points on $l_{nqx}$ must have their $\gamma$-balls
infected by the time an outburst whose radius exceeds $\gamma$ has
occurred in each of the areas $A_k$ $(k=1,\ldots ,12)$. Hence
$B(z,\gamma)\subset S_{t_0+\max\{E_1,\ldots ,E_{12}\}}$ for all
$z\in l_{nqx}$ and we have proved that $\tilde{T}(l_{nqx})-t_0\leq
\max\{E_1,\ldots ,E_{12}\}$ for all times $t_0$ such that there
exists a point on $l_{nqx}$ whose entire $\gamma$-ball is infected
at time $t_0$. The above reasoning easily generalizes to $d\geq
3$. We obtain

\begin{equation} \label{eq:allaeb}
\tilde{T}(l_{nqx})-t_0\leq \max\{E_1,\ldots ,E_{k_0}\},
\end{equation}

\noindent where $k_0=k_0(d)$ and $E_k\sim
\textrm{Exp}(\lambda(d))$.\medskip

\noindent Now, to prove the first inequality in
(\ref{eq:lstildeTtxmmuli}) let $q\in\mathbb{Q}$ be small enough to
ensure (\ref{eq:allaeb}) and let $\tilde{\varphi}_{nq}(x)$ be the
time when the $\gamma$-balls of all points on $l_{nqx}$ are
infected counting from the time when the $\gamma$-ball around
$nqx$ is infected, that is, $\tilde{\varphi}_{nq}(x)=
\tilde{T}(l_{nqx})-\tilde{T}(nqx)$. By (\ref{eq:allaeb}) there
exist i.i.d.\ random variables $E_1,\ldots ,E_{k_0}$ such that

\begin{equation} \label{eq:varphileqrn}
\tilde{\varphi}_{nq}(x)\leq \max\{E_1,\ldots ,E_{k_0}\}:=M_{k_0}.
\end{equation}

\noindent Since $M_{k_0}\leq \sum_{k=1}^{k_0}E_k$, it follows that
$\textrm{E}[M_{k_0}]\leq k_0\textrm{E}[E_1]<\infty$ and, hence,
for all $\varepsilon>0$,
$$
\sum_{n=0}^{\infty}P(M_{k_0}>n\varepsilon)<\infty.
$$
\noindent Thus, by (\ref{eq:varphileqrn}) and the Borel-Cantelli
lemma,

\begin{equation} \label{eq:phitildenoll}
\lim_{n\rightarrow \infty} \frac{\tilde{\varphi}_{nq}(x)}{n}=0
\quad \textrm{a.s.}
\end{equation}

\noindent Let $n_t$ be such that $t\in[n_tq, (n_t+1)q)$. Clearly
$\tilde{T}(tx)\leq \tilde{T}(n_tqx)+ \tilde{\varphi}_{nq}(x)$.
Hence

\begin{eqnarray}
\limsup_{t\rightarrow \infty}\frac{\tilde{T}(tx)}{t} & \leq &
\limsup_{t\rightarrow \infty} \frac{\tilde{T}(n_tqx)+
\tilde{\varphi}_{n_tq}(x)}{n_tq}
\nonumber \\
 & = & \lim_{n\rightarrow \infty} \frac{\tilde{T}(nqx)+\tilde{\varphi}_{nq}(x)}{nq}
  \label{eq:kontgv1} \\[0.3cm]
 & = & \mu(x) \nonumber,
\end{eqnarray}

\noindent where the last equality follows from
(\ref{eq:tildeTncx}) and (\ref{eq:phitildenoll}).\medskip

\noindent To prove the last inequality in
(\ref{eq:lstildeTtxmmuli}), let $z_{nqx}$ be the first point on
$l_{nqx}$ whose $\gamma$-ball is infected and let
$\tilde{\psi}_{nq}(x)$ be the time from when this occurs until the
infection has invaded the $\gamma$-balls of all points on
$l_{nqx}$, that is, $\tilde{\psi}_{nq}(x)=\tilde{T}(l_{nqx})
-\tilde{T}(z_{nqx})$. Using (\ref{eq:allaeb}) and Borel-Cantelli
it can be seen that

\begin{equation} \label{eq:psitildenoll}
\lim_{n\rightarrow \infty} \frac{\tilde{\psi}_{nc}(x)}{n}=0 \quad
\textrm{a.s.}
\end{equation}

\noindent Since $\tilde{T}(tx)+ \tilde{\psi}_{n_tq}(x)\geq
\tilde{T}(n_tqx)$, we have

\begin{eqnarray}
\liminf_{t\rightarrow \infty} \frac{\tilde{T}(tx)}{t}& \geq &
\liminf_{t\rightarrow \infty}
\frac{\tilde{T}(n_tqx)-\tilde{\psi}_{n_tq}(x)}{(n_t+1)}
\nonumber \\
 & = & \lim_{n\rightarrow \infty} \frac{\tilde{T}(nqx)-\tilde{\psi}_{nq}(x)}{(n+1)q}
 \label{eq:kontgv2} \\[0.3cm]
 & = & \mu(x) \nonumber,
\end{eqnarray}

\noindent where the last equality is a consequence of
(\ref{eq:tildeTncx}) and (\ref{eq:psitildenoll}). Thereby
(\ref{eq:lstildeTtxmmuli}) is established and the proposition
follows. $\hfill \Box$\bigskip

\noindent It follows from Proposition \ref{prop:kontpar} that
$\mu(cx)=c\mu(x)$. This implies that $\mu(x)=|x|\mu(\hat{x})$,
where $\hat{x}=x/|x|$, the unit vector in direction $x$. Due to
rotational invariance of $\mathbb{R}^d$ and the model it is clear
that $\mu(\hat{x})=\mu(\hat{y})$ for all $x,y\in \mathbb{R}^d$.
Thus we can define a constant
$$
\mu:=\mu((1,0,\ldots ,0))
$$
representing the asymptotic time it takes for the infection to
travel a unit vector in an arbitrary direction. By Proposition
\ref{prop:0muinf}, $\mu \in(0,\infty)$. We end up with the simple
relation $ \mu(x)=|x|\mu$ valid for all $x\in
\mathbb{R}^d$.\medskip

\noindent To summarize the results obtained in the present
section, we have deduced that there is a real number $\mu>0$ such
that, for each fixed $x\in \mathbb{R}^d$, almost surely

\begin{equation} \label{eq:step1}
\lim_{t\rightarrow \infty}\frac{\tilde{T}(tx)}{t}=|x|\mu.
\end{equation}

\section{Proof of the shape theorem}

\noindent The shape theorem asserts that $S_t\approx
tB(0,\mu^{-1})$ for large $t$ and in view of (\ref{eq:step1}) this
is indeed what to expect. However, it remains to show that the
linear growth stipulated in (\ref{eq:step1}) holds for all
directions simultaneously. To this end we will need the following
lemma, which asserts that with high probability the infected
region in a process emanating from a point $y$ will eventually
contain a ball centered at $y$ with radius proportional to time.

\begin{lemma} \label{lemma:Stilb}
For small $\delta>0$, there is a constant $c\in (0,\infty)$ and a
time $s_0$ such that, for any $y\in \mathbb{R}^d$ and $s'\geq 0$,
we have
$$
P(B(y,s\delta)\not\subset S_{s'+s}^{(y,s')})\leq e^{-cs}
$$
\noindent if $s>s_0$.
\end{lemma}

\noindent \emph{Proof of Lemma \ref{lemma:Stilb}:} Due to shift
invariance of the model it suffices to prove the lemma for $y=0$
and $s'=0$, that is, it suffices to show that

\begin{equation} \label{eq:BtdinSt}
P(B(0,s\delta)\not\subset S_s)\leq e^{-cs}
\end{equation}

\noindent for large $s$. To achieve this, fix a real number
$\alpha>0$ and partition $\mathbb{R}^d$ into cubes centered at the
points $\alpha \mathbb{Z}^d$ and with vertices $(\alpha/2,\ldots
,\alpha/2)+\alpha \mathbb{Z}^d$. Furthermore, let
$\hat{T}_\alpha(x)$, $x\in \alpha \mathbb{Z}^d$, denote the time
when the entire cube centered at $x$ is infected. For small
$\alpha$ and $\delta$ we will show that there is a positive
constant $c'$ such that the estimate

\begin{equation} \label{eq:PThatbeg}
P(\hat{T}_\alpha(x)>a|x|/\delta)\leq e^{-c'a|x|}
\end{equation}

\noindent holds simultaneously for all $a\geq 1/d$ and all
$x\in\mathbb{Z}^d$. Given this estimate the lemma is readily
established: Let $C(s\delta+\alpha)$ denote the cube with side
length $2(s\delta+\alpha)$ centered at the origin and write
$C_{\alpha}(s\delta+\alpha)=\alpha \mathbb{Z}^d \cap
C(s\delta+\alpha)$. Note that $B(0,s\delta)$ is contained in the
union of all $\alpha$-cubes whose centers are in
$C_{\alpha}(s\delta+\alpha)$. Thus, if $B(0,s\delta)\not\subset
S_s$, then not all $\alpha$-cubes with centers in
$C_{\alpha}(s\delta+\alpha)$ are infected at time $s$, that is,

\begin{eqnarray*}
P(B(0,s\delta)\not\subset S_s) & \leq & P\left(\bigcup_{x\in
C_{\alpha}(s\delta+\alpha)}\{\hat{T}_\alpha(x) >s\}
\right) \\
 & \leq & \sum_{x\in C_{\alpha}(s\delta+\alpha)}
 P(\hat{T}_\alpha(x)>s).
\end{eqnarray*}

\noindent Trivially,
$$
P(\hat{T}_\alpha(x)>s)=P\left(\hat{T}_\alpha(x)>
\frac{s\delta}{|x|}\cdot \frac{|x|}{\delta}\right).
$$
\noindent For $x\in C_{\alpha}(s\delta+\alpha)$, we have $|x|\leq
\sqrt{d}(s\delta+\alpha)$ and, since $s\delta/\sqrt{d}
(s\delta+\alpha)\rightarrow 1/\sqrt{d}$ as $s\rightarrow \infty$,
it holds that $s\delta/|x|\geq 1/d$ for large $s$. Hence, if $s$
is large and $x\in C_{\alpha}(s\delta+\alpha)$ it follows from
(\ref{eq:PThatbeg}) that
$$
P\left(\hat{T}_\alpha(x)> \frac{s\delta}{|x|}\cdot
\frac{|x|}{\delta}\right)\leq e^{-c's\delta}
$$
\noindent and, consequently, for large $s$,
$$
P(B(0,s\delta)\not\subset S_s)\leq \sum_{x\in
C_{\alpha}(s\delta+\alpha)} e^{-c's\delta}.
$$
\noindent Since the number of points in
$C_{\alpha}(s\delta+\alpha)$ grows only polynomially in $s$ there
is a time $s_0$ such that, for $s>s_0$, we have
$P(B(0,s\delta)\not\subset S_s)\leq e^{-sc'\delta/2}$ as
desired.\medskip

\noindent It remains to prove (\ref{eq:PThatbeg}). Fix $x\in
\alpha \mathbb{Z}^d$. By Lemma \ref{lemma:txyb} and Remark 3.1, if
$\alpha$ is small, say $\alpha \leq \gamma/10$, we have
$$
\hat{T}_\alpha(x)\leq \sum_{k=1}^{2\lceil|x|\gamma^{-1}\rceil}E_k,
$$
\noindent where $\{E_k\}$ are i.i.d.\ exponential variables with
parameter $\lambda$. Thus, it suffices to find $c'>0$ such that

\begin{equation} \label{eq:anka}
P\left(\sum_{k=1}^{2\lceil|x|\gamma^{-1}\rceil}E_k>\frac{a|x|}{\delta}\right)\leq
e^{-c'a|x|}.
\end{equation}

\noindent To this end, write $\lceil|x|\gamma^{-1}\rceil=m$ and
introduce the notation $\Gamma_{2m}:=\sum_{k=1}^{2m}E_k$. Using
Markov's inequality and the fact that $\Gamma_{2m}\sim
\textrm{Gamma}(2m,\lambda)$, we obtain
$$
e^{\theta
a|x|\delta^{-1}}P\left(\Gamma_{2m}>a|x|\delta^{-1}\right)\leq
\textrm{E}[e^{\theta \Gamma_{2m}}]=(1-\lambda \theta)^{-2m}
$$
\noindent for $\theta \in(0,\lambda^{-1})$. Thus

\begin{equation} \label{eq:tuppen}
P\left(\Gamma_{2m}>a|x|\delta^{-1}\right)\leq
\exp\left\{-a|x|\left(\theta
\delta^{-1}+\frac{2m}{a|x|}\log(1-\theta \lambda)\right)\right\}.
\end{equation}

\noindent We may assume that $|x|>\gamma -\alpha$ since, for $x\in
B(0,\gamma -\alpha)$, the $\alpha$-cube centered at $x$ is
contained in $B(0,\gamma)$, implying that the left hand side in
(\ref{eq:PThatbeg}) equals zero and hence (\ref{eq:PThatbeg}) is
trivially true in this case. For $|x|>\gamma -\alpha$, the
quotient $m/|x|=\lceil|x|\gamma^{-1}\rceil/|x|$ is bounded by
$2\gamma^{-1}$. Substituting this in (\ref{eq:tuppen}) and also
using the fact that $a\geq 1/d$ yields
$$
P\left(\Gamma_{2m}>a|x|\delta^{-1}\right)\leq
e^{-a|x|f_{\delta}(\theta)},
$$
\noindent where $f_{\delta}(\theta)=\theta
\delta^{-1}+4d\gamma^{-1}\log(1-\lambda \theta)$. Now,
$f_{\delta}(0)=0$ and
$f_{\delta}'(0)=\delta^{-1}-4d\gamma^{-1}\lambda$. Thus, if
$\delta$ is so small that $f_{\delta}'(0)>0$, then we can pick
$\theta$ small and get $f_{\delta}(\theta)>0$. This proves
(\ref{eq:anka}). $\hfill \Box$\bigskip

\noindent Finally, equipped with the above lemma and the results
from Section 4, we are ready to prove the shape theorem.\bigskip

\noindent \emph{Proof of Theorem \ref{th:aform}:} Fix
$\varepsilon\in (0,\mu^{-1})$. We will prove the theorem in two
steps:

\begin{itemize}
\item[(i)] There is almost surely a time $T^1$ such that
$(1-\varepsilon)tB(0,\mu^{-1})\subset S_t$ for $t>T^1$.
\item[(ii)] There is almost surely a time $T^2$ such that
$S_t\subset(1+\varepsilon)tB(0,\mu^{-1})$ for $t>T^2$.
\end{itemize}

\noindent As for (i) we will, in fact, prove the following:

\begin{itemize}
\item[(i')] There is almost surely a time $\acute{T}^1$ such that
$(1-\varepsilon/2)tB(0,\mu^{-1})\subset S_t$ for $t>\acute{T}^1$,
$t\in \mathbb{N}$.
\end{itemize}

\noindent From (i') it follows that $(1-\varepsilon/2)\lfloor
t\rfloor B(0,\mu^{-1})\subset S_{\lfloor t\rfloor}$ for
$t>\acute{T}^1+1$, where $\lfloor \cdot\rfloor$ is the integer
part function. Since $S_{\lfloor t\rfloor}\subset S_t$ for all
$t$, we obtain $(1-\varepsilon/2)\lfloor t\rfloor
B(0,\mu^{-1})\subset S_t$. But, for large $t$, the ball with
radius $(1-\varepsilon/2)\lfloor t\rfloor \mu^{-1}$ contains the
ball with radius $(1-\varepsilon)t\mu^{-1}$, and hence (i) follows
from (i').\medskip

\noindent To prove (i'), note that, since
$(1-\varepsilon/2)B(0,\mu^{-1})$ is compact, there are points
$x_1,\ldots ,x_n\in (1-\varepsilon/2)B(0,\mu^{-1})$ such that
$$
(1-\varepsilon/2)B(0,\mu^{-1})\subset \bigcup_{i=1}^n B(x_i,\delta
\varepsilon /4),
$$
\noindent where $\delta>0$ is chosen small enough to ensure that
Lemma \ref{lemma:Stilb} holds. Clearly

\begin{equation} \label{eq:(i)3}
(1-\varepsilon/2)tB(0,\mu^{-1})\subset \bigcup_{i=1}^n
B(tx_i,t\delta \varepsilon /4).
\end{equation}

\noindent Since $|x_i|\leq (1-\varepsilon/2)\mu^{-1}$ it follows
from (\ref{eq:step1}) that, for each $i$, $\lim_{t\rightarrow
\infty}\tilde{T}(tx_i)/t\leq 1-\varepsilon/2$ almost surely. This
implies that almost surely $\tilde{T}(tx_i)\leq
t(1-\varepsilon/4)$ for each $i$ if $t$ is large, that is,

\begin{equation} \label{eq:(i)1}
B(tx_i,\gamma)\subset S_{t(1-\varepsilon/4)}
\end{equation}

\noindent for large $t$. Furthermore, by Lemma \ref{lemma:Stilb},
$$
\sum_{t\in \mathbb{N}}P\left(B(tx_i,t\delta
\varepsilon/4)\not\subset
S_t^{(tx_i,t(1-\varepsilon/4))}\right)\sim \sum_{t\in \mathbb{N}}
e^{-ct\varepsilon/4}<\infty.
$$
\noindent Thus, by the Borel-Cantelli lemma, for large integer
times we have almost surely

\begin{equation} \label{eq:(i)2}
B(tx_i,t\delta \varepsilon/4)\subset
S_t^{(tx_i,t(1-\varepsilon/4))}.
\end{equation}

\noindent Now, for each $i$, let $T(i)$ be such that both
(\ref{eq:(i)1}) and (\ref{eq:(i)2}) hold for $t>T(i)$, $t\in
\mathbb{N}$, and define $\acute{T}^1=\max\{T(i)\}$. For
$t>\acute{T}^1$, we have $B(tx_i,t\delta \varepsilon/4)\subset
S_t$ for all $i$ and $t\in \mathbb{N}$. Using (\ref{eq:(i)3}) this
implies that $(1-\varepsilon/2)tB(0,\mu^{-1})\subset S_t$ for all
integer times larger than $T_1'$, as desired.\medskip

\noindent Moving on to (ii), let $A$ be an annulus of width
$\mu^{-1}\varepsilon/2$ surrounding $(1+\varepsilon/2)\linebreak
B(0,\mu^{-1})$, i.e.
$$
A=(1+\varepsilon) B(0,\mu^{-1})\backslash (1+\varepsilon/2)
B(0,\mu^{-1}),
$$
\noindent and let $\tilde{S}_t$ be the set of points whose entire
$\gamma$-ball is infected at time $t$, i.e.
$$
\tilde{S}_t=\{x;\hspace{0.1cm} B(x,\gamma)\subset S_t\}.
$$
\noindent Points in $\tilde{S}_t$ will be referred to as strongly
infected at time $t$. We begin by showing that almost surely
$\tilde{S}_t\cap tA=\emptyset$ for large $t$.\medskip

\noindent Let $\delta$ be small enough to ensure that Lemma
\ref{lemma:Stilb} holds and pick $x_1,\ldots ,x_n\in A$ such that

\begin{equation} \label{eq:Rt}
A\subset \bigcup_{i=1}^{n}B(x_i,\delta\varepsilon/8).
\end{equation}

\noindent By (\ref{eq:step1}), $\lim_{t\rightarrow
\infty}\tilde{T}(tx_i)/t\geq 1+\varepsilon/2$ for each $i$. Hence,
for each $i$, almost surely $\tilde{T}(tx_i)\geq
t(1+\varepsilon/4)$ for large $t$, implying that, for every $c\in
(0,1)$, we have

\begin{equation} \label{eq:(ii)1}
P\left(B(tx_i,\gamma)\subset S_{t(1+\varepsilon/4)}\textrm{ for
some }i=1,\ldots ,n \right)\leq c
\end{equation}

\noindent if $t$ is large. The idea of the proof is that, if $tA$
contains strongly infected points for large $t$, then with high
probability some point $tx_i$ will be strongly infected within
time $t\varepsilon/4$ and this conflicts with (\ref{eq:(ii)1}). To
formalize this intuition, let
$$
p=P\left(\tilde{S}_t\cap tA\neq \emptyset \textrm{ for arbitrarily
large }t\right)
$$
\noindent and assume for contradiction that $p>0$. For fixed $t$,
write $E_t=\{t'\geq t; \hspace{0.1cm} \tilde{S}_{t'}\cap t'A\neq
\emptyset\}$ and define
$$
T_t=\left\{ \begin{array}{lc}
            \inf E_t & \textrm{if } E_t\neq \emptyset,\\
            \infty & \textrm{if } E_t=\emptyset.
            \end{array} \right.
$$
\noindent Note that $P(T_t<\infty)\geq p>0$ for each $t$.
Consequently, we can condition on the event that $T_t<\infty$ and
pick $y_t$ uniformly on $\tilde{S}_{T_t}\cap T_tA$. By Lemma
\ref{lemma:Stilb},
$$
P\left(B(y_t,T_t\delta\varepsilon/4)\subset
S_{T_t(1+\varepsilon/4)}^{(y_t,T_t)}\Big|\hspace{0.1cm} T_t<\infty
\right)\geq 1-e^{-ct\varepsilon/4},
$$
\noindent that is, the ball with radius $T_t\delta \varepsilon/4$
around $y_t$ is with high probability infected at time
$T_t(1+\varepsilon/4)$ in a process started from $y_t$ at time
$T_t$. If $t$ is large this implies that the ball with radius
$T_t\delta \varepsilon/8$ is strongly infected with high
probability. Formally,
$$
P\left(B(y_t,T_t\delta\varepsilon/8)\subset
\tilde{S}_{T_t(1+\varepsilon/4)}^{(y_t,T_t)}\Big|\hspace{0.1cm}
T_t<\infty \right)\geq 1-e^{-ct\varepsilon/4}.
$$
\noindent Since $y_t$ is strongly infected at time $T_t$ we obtain
$$
P\left(B(y_t,T_t\delta\varepsilon/8)\subset
\tilde{S}_{T_t(1+\varepsilon/4)}\Big|\hspace{0.1cm} T_t<\infty
\right)\rightarrow 1 \quad \textrm{as } t\rightarrow \infty.
$$
\noindent Now, by (\ref{eq:Rt}), $T_tA$ is covered by the balls
$B(T_tx_i,T_t\delta\varepsilon/8)$. Hence, since $y_t\in T_tA$, we
can find at least one point $x_i$ such that $T_tx_i\in
B(y_t,T_t\delta\varepsilon/8)$ and, consequently,
$$
P\left(B(T_tx_i,\gamma)\subset S_{T_t(1+\varepsilon/4)}\textrm{
for some }i=1,\ldots ,n|\hspace{0.1cm}
T_t<\infty\right)\rightarrow 1 \quad \textrm{as } t\rightarrow
\infty.
$$
\noindent Pick $t$ large so that the above probability is greater
than 1/2 and so that (\ref{eq:(ii)1}) holds for $c=p/4$. Use the
fact that $P(T_t<\infty)\geq p$ for each $t$ to obtain
$$
\begin{array}{ll}
 & P\left(\exists t'\geq t \textrm{ such that }
 B(t'x_i,\gamma)\subset S_{t'(1+\varepsilon/4)} \textrm{ for some }i=1,
 \ldots ,n\right)\\
 [0.2cm] \geq & P\left(B(T_tx_i,\gamma)\subset S_{T_t(1+\varepsilon/4)}\textrm{ for some
}i=1,\ldots ,n|\hspace{0.1cm} T_t<\infty\right)
P\left(T_t<\infty\right)\\[0.2cm]
> & p/2.
\end{array}
$$

\noindent This contradicts (\ref{eq:(ii)1}). Hence we must have
$p=0$, that is, almost surely $\tilde{S}_t\cap tA=\emptyset$ for
large $t$. It remains to show that (ii) follows from this. To this
end, let $\Gamma$ denote the set of outbursts that occur in $tA$
for some $t$ during the progress of the growth, that is, $\Gamma$
is the set of outbursts $(X_n,T_n)$ such that $X_n\in T_nA$.
Furthermore, let $\Gamma_{\gamma}$ be those outbursts in $\Gamma$
whose radius is at least $\gamma$.  Assume that (ii) fails so that
with positive probability
$S_t\cap[(1+\varepsilon)tB(0,\mu^{-1})]^c\neq \emptyset$ for
arbitrarily large $t$. This implies that $P(|\Gamma|=\infty)>0$.
Now, each time an outburst takes place we condition on the process
up to that time, including the fact that an outburst takes place
at that time at that location, but excluding the radius of the
outburst. Since the radii of the outbursts are i.i.d. it follows
from Levy's version of the Borel-Cantelli lemma (see Williams
(1991), section 12.5) that we can not have $|\Gamma|=\infty$ and
$|\Gamma_\gamma|<\infty$. Hence the fact that
$P(|\Gamma|=\infty)>0$ implies that
$P(|\Gamma_{\gamma}|=\infty)>0$ as well. But, if
$|\Gamma_{\gamma}|=\infty$, then the region $tA$ must contain
strongly infected points for arbitrarily large $t$, and hence we
have derived a contradiction.\medskip

\noindent At this point, (i) and (ii) are established and all that
remains is to note that, for $t>\max\{T^1,T^2\}$, we have
$$
(1-\varepsilon)tB(0,\mu^{-1})\subset S_t\subset
(1+\varepsilon)tB(0,\mu^{-1}).
$$
\noindent The shape theorem is proved. $\hfill \Box$\bigskip

\noindent We end by proving Theorem \ref{th:oberS0}, that is, by
proving that the asymptotic shape is independent of the state at
time zero.\bigskip

\noindent \emph{Proof of Theorem \ref{th:oberS0}:} For an
arbitrary set $\Pi\subset \mathbb{R}^d$, introduce generalized
versions of the quantities $S_t$ and $S_t^{(x,s)}$ as follows: Let
$S_t^\Pi$ denote the region infected at time $t$ in a process
started from the set $\Pi$ and let $S_t^{(\Pi,s)}$ denote the
region infected at time $t$, $t\geq s$, in a process started from
the set $\Pi$ at time $s$. We will show that, for any
$\varepsilon>0$, almost surely

\begin{equation}\label{eq:smomrinkl}
(1-\varepsilon)tB(0,\mu^{-1})\subset \frac{S_t^\Pi}{t}\subset
(1+\varepsilon)tB(0,\mu^{-1})
\end{equation}

\noindent if $t$ is large. To this end, let $\tau$ be the first
time when $\Pi$ is contained in $S_t$, and let $\tau'$ be the
first time counting from $\tau$ when $B(0,\gamma)$ is contained in
$S^{(\Pi,\tau)}_t$, that is,
$$
\tau=\inf\{t;\hspace{0.1cm} \Pi\subset S_t\}
$$
\noindent and
$$
\tau'=\inf\{t;\hspace{0.1cm} B(0,\gamma)\subset
S^{(\Pi,\tau)}_{\tau+t}\}.
$$
\noindent If $\Pi$ is bounded and has strictly positive Lebesgue
measure these times are both finite with probability 1.
Furthermore,
$$
S^{(B(0,\gamma),\tau+\tau')}_{\tau+\tau'+t}\subset
S^{(\Pi,\tau)}_{\tau+\tau'+t}\subset S_{\tau+\tau'+t}
$$
\noindent which implies that
$$
\frac{S^{(B(0,\gamma),\tau+\tau')}_{\tau+\tau'+t}}{t}\subset
\frac{\tau'+t}{t}\cdot
\frac{S^{(\Pi,\tau)}_{\tau+\tau'+t}}{\tau'+t}\subset
\frac{\tau+\tau'+t}{t}\cdot \frac{S_{\tau+\tau'+t}}{\tau+\tau'+t}.
$$
\noindent The quotients $(\tau'+t)/t$ and $(\tau+\tau'+t)/t$ tend
to 1 as $t\rightarrow \infty$ and by Theorem \ref{th:aform} the
asymptotic shape for the processes
$S^{(B(0,\gamma),\tau+\tau')}_{\tau+\tau'+t}/t$ and
$S_{\tau+\tau'+t}/(\tau+\tau'+t)$ is $B(0,\mu^{-1})$. This
establishes (\ref{eq:smomrinkl}). $\hfill \Box$

\section{Acknowledgement}

\noindent I thank my thesis advisor Olle H\"{a}ggstr\"{o}m for
proposing the problem, for superb guidance during the progress of
the work and for critical readings of the manuscript. I also thank
Anders Martin-L\"{o}f for many stimulating discussions. Finally, I
tnak an anonymous referee for his or her careful reading and
constructive comments.

\section*{References}

\noindent Boivin, D. (1990): First passage percolation: the
stationary case, \emph{Prob. Th. Rel. Fields} \textbf{86},
491-499.\medskip

\noindent Bramson, M. and Griffeath, D. (1981): On the
Williams-Bjerknes tumour growth model I, \emph{Ann. Prob.}
\textbf{9}, 173-185.\medskip

\noindent Cox, J.T. and Durrett, R. (1981): Some limit theorems
for percolation processes with necessary and sufficient
conditions, \emph{Ann. Prob.} \textbf{9}, 583-603.\medskip

\noindent Durrett, R. (1988): \emph{Lecture Notes on Particle
Systems and Percolation}, Wadsworth $\&$ Brooks/Cole.\medskip

\noindent Eden, M. (1961): A two-dimensional growth process,
\emph{Proceedings of the 4th Berkeley symposium on mathematical
statistics and probability} vol. \textbf{IV}, 223-239, University
of California press.\medskip

\noindent Howard C.D. and Newman C.M. (1997): Euclidean models of
first-passage percolation, \emph{Prob. Th. Rel. Fields}
\textbf{108}, 153-170.\medskip

\noindent Kesten, H. (1986): Aspects of first-passage percolation,
\emph{Lecture Notes in Mathematics}, vol. \textbf{1180}, 125-264,
Springer .\medskip

\noindent Lee, T. and Cowan, R. (1994): A stochastic tessellation
of digital space, \emph{Mathematical morphology and its
applications to image processing}, 217-224, Kluwer.\medskip

\noindent Liggett, T.M. (1985): An improved subadditive ergodic
theorem, \emph{Ann. Prob.} \textbf{13}, 1279-1285.\medskip

\noindent Penrose, M. (2001): Random parking, sequential
adsorption and the jamming limit, \emph{Commun. Math. Phys.}
\textbf{218}, 153-176.\medskip

\noindent Richardson, D. (1973): Random growth in a tessellation,
\emph{Proc. Cambridge Phil. Soc.} \textbf{74}, 515-528.\medskip

\noindent Williams, D. (1991): \emph{Probability with
Martingales}, Cambridge University Press.\medskip

\noindent Williams, T. and Bjerknes R. (1972): Stochastic model
for abnormal clone spread through epithelial basal layer,
\emph{Nature} \textbf{236}, 19-21.

\end{document}